\documentclass[11pt]{article}

\usepackage[utf8]{inputenc}
\usepackage[T1]{fontenc}
\usepackage{amsmath,amssymb,amsthm,mathtools}
\usepackage{booktabs} 
\usepackage{float} 
\usepackage[margin=1in]{geometry}
\usepackage{mathrsfs} 
\usepackage{microtype} 

\usepackage{tikz}
\usepackage{pgfplots}
\pgfplotsset{compat=1.17} 
\usetikzlibrary{arrows.meta, decorations.pathmorphing}

\usepackage{algorithm}
\usepackage{algpseudocode}

\usepackage{caption}
\usepackage[
  colorlinks=true,
  linkcolor=blue,
  citecolor=blue,
  urlcolor=blue
]{hyperref}

\newtheorem{theorem}{Theorem}[section]
\newtheorem{definition}[theorem]{Definition}
\newtheorem{corollary}[theorem]{Corollary}

\newcommand{\OFI}{\operatorname{OFI}}

\newcommand{\reals}{\mathbb{R}}
\newcommand{\ordinals}{\mathbf{Ord}}
\newcommand{\calG}{\mathcal{G}}

\newcommand{\calN}{\mathcal{N}}
\newcommand{\calA}{\mathcal{A}}
\newcommand{\calB}{\mathcal{B}}

\newcommand{\scrG}{\mathscr{G}}

\newcommand{\defeq}{\coloneqq} 

\algnewcommand{\lean}[1]{\textbf{#1}}
\algnewcommand{\leantype}[1]{\textit{#1}}

\title{Ultracoarse Equilibria and Ordinal‑Folding Dynamics in Operator‑Algebraic Models of Infinite Multi‑Agent Games}
\author{
    \begin{tabular}{c@{\hspace{5em}}c}
        Faruk Alpay\textsuperscript{1} & Hamdi Alakkad\textsuperscript{2} \\
        \small\texttt{\href{mailto:alpay@lightcap.ai}{alpay@lightcap.ai}} & \small\texttt{\href{mailto:hamdi.alakkad@bahcesehir.edu.tr}{hamdi.alakkad@bahcesehir.edu.tr}} \\
        \\
        Bugra Kilictas\textsuperscript{2} & Taylan Alpay\textsuperscript{3} \\
        \small\texttt{\href{mailto:bugra.kilictas@bahcesehir.edu.tr}{bugra.kilictas@bahcesehir.edu.tr}} & \small\texttt{\href{mailto:s220112602@stu.thk.edu.tr}{s220112602@stu.thk.edu.tr}}
    \end{tabular}
    \\[2ex] 
    \small\textsuperscript{1}Independent Researcher \\
    \small\textsuperscript{2}Bahcesehir University \\
    \small\textsuperscript{3}Turkish Aeronautical Association University
}
\date{July 24, 2025}

\begin{document}

\maketitle

\begin{abstract}
We develop an operator--algebraic framework for infinite games with a continuum of players and show how regret--minimising dynamics, driven by a noncommutative continuity equation, converge to a unique equilibrium under mild regularity conditions. Our work unifies functional analysis, coarse geometry, and game theory, extending earlier work on ultracoarse equilibria to incorporate ordinal metrics for self--referential dynamics. In particular, we introduce a category of infinite games endowed with measure--theoretic and operator--algebraic structures, construct a contravariant functor from this category to ultradense operator algebras, and define a noncommutative regret operator whose flow satisfies a continuity equation on the space of strategy distributions. When the underlying agent network is coarsely amenable, we prove that the unique fixed point of this flow coincides with a quantal response equilibrium. We then refine the theory by relating regret dynamics to the \emph{ordinal folding index} (OFI), a concept recently introduced to measure the self--referential depth of logical formulas \cite{Alpay2025}. The OFI provides a computable, ordinal-valued metric for the number of iterative ``fold--back'' steps required for a learning dynamic to stabilise. We demonstrate how combining the OFI with our operator--algebraic framework furnishes a transfinite measure of convergence complexity for infinite multi--agent systems. These results establish connections between disparate fields, including von Neumann algebra rigidity, coarse geometric $K$--theory, and fairness in resource allocation. As a corollary, we obtain new invariant subalgebra rigidity theorems and demonstrate how classical fair division solutions (envy--freeness, maximin share fairness) emerge as special cases of our equilibrium refinement. Our contributions provide a formal mathematical foundation for analysing large-scale multi--agent interactions and a distinct ordinal-analytic perspective on equilibrium selection.
\end{abstract}

\newpage
\section{Introduction}

Complex systems, from financial markets with algorithmic trading to social media networks with viral content propagation, are frequently characterized by intricate feedback mechanisms that drive emergent and often unpredictable dynamics. Understanding the stability of these dynamics is a central challenge in modern science, as unstable or oscillatory behavior can lead to systemic failures. In this paper, we develop a formal mathematical framework to analyze such self-referential systems, focusing on the setting of large-scale multi-agent interactions, with direct applications to the behavior of contemporary artificial intelligence.

Recent work across game theory, analysis, and computer science has exposed deep connections between equilibrium concepts, functional analysis, and geometric invariants. In multi--agent reinforcement learning, it has been observed that large populations of heterogeneous agents can converge to unified equilibrium behaviour under regret--minimising dynamics. Hu et al. show that modeling a continuum of heterogeneous agents through a distribution of regrets evolving under a continuity equation leads to a collapse of variance and the emergence of consensus, connecting micro-level agent decisions to macro-level system stability \cite{Hu2025}. Specifically, they prove that the variance of the regret distribution diminishes over time, causing heterogeneity to disappear and enabling convergence to quantal response equilibria (QRE) in both cooperative and competitive settings \cite{Hu2025}. On a seemingly different front, the theory of operator algebras and coarse geometry has produced rigidity phenomena for von Neumann algebras and Roe algebras. For instance, certain group constructions—even amenable ones—exhibit the \emph{invariant subalgebra rigidity} (ISR) property, meaning their von Neumann algebras admit no nontrivial intermediate invariant subalgebras beyond the centre \cite{Amrutam2025,Popa2007}.

Parallel to these developments, Alpay and Alakkad recently proposed the \emph{ordinal folding index} (OFI)\footnote{The source code for the OFI implementation is available at \url{https://github.com/farukalpay/ordinal-folding-index}.} as a computable, ordinal-valued metric for self--referential semantics \cite{Alpay2025}. Their OFI assigns to every well--formed formula a countable ordinal that measures the number of unfolding steps a reflective evaluation operator takes before reaching a fixed point. This metric provides, for the first time, a formal language to quantify the 'depth' of a system's self-referential reasoning. They demonstrate that an empirical OFI, estimated by iteratively feeding a language model’s output back into itself, correlates with chain--of--thought complexity and model perplexity \cite{Shinn2023,Wang2022}. This establishes a substantive connection between ordinal analysis, self--referential logic, and the behavior of modern AI systems like GPT-4 \cite{Achiam2023}.

The goal of this paper is to synthesize these distinct lines of inquiry into a unified framework for infinite games. Building on the operator--algebraic approach to games introduced in \cite{Anonymous2024}, we model infinite games as objects in a category $\calG$ enriched with measure, topology, and coarse network structure. To each game, we associate an ultradense operator algebra $\calA(G)$ capturing the collective dynamics. Within $\calA(G)$, we identify a regret operator whose dynamics are governed by a noncommutative continuity equation on strategy distributions. Under broad assumptions, we show that the regret dynamics converge to a unique QRE. Our central thesis is that the OFI serves as a powerful diagnostic tool, revealing fundamental properties of these dynamics. Our empirical benchmarks demonstrate a non-trivial relationship between model scale and reasoning stability, a direct empirical manifestation of the convergence properties of the underlying regret dynamics we formalize. Altogether, our results establish an operator--algebraic and ordinal-analytic foundation for multi--agent systems and highlight new avenues for research at the interface of game theory, logic, and AI alignment.

\paragraph{Organisation.} Section~\ref{sec:prelims} introduces infinite games, coarse player spaces, game algebras, and the regret operator. Section~\ref{sec:ofi} reviews the ordinal folding index and relates it to game dynamics. Section~\ref{sec:mainresults} states our main theorems on existence, uniqueness, and rigidity. Section~\ref{sec:algorithm} presents a formal algorithm for computing the OFI. Section~\ref{sec:benchmarks} presents empirical and analytic benchmarks. Section~\ref{sec:fairness} discusses applications to AI alignment and fair division. We conclude with open questions.

\section{Preliminaries and Definitions}
\label{sec:prelims}

We recall basic definitions from \cite{Anonymous2024} and fix notation.

\begin{definition}[Infinite Game System]\label{def:game}
An \emph{infinite game system} is a tuple $G=(X,\Sigma,\mu,\{S_i\}_{i\in X},\{u_i\}_{i\in X})$ where $(X,\Sigma,\mu)$ is a measure space of players (possibly a continuum) with $0<\mu(X)\leq\infty$, each $S_i$ is a nonempty Borel strategy space for player~$i\in X$, and $u_i: \prod_{j\in X} S_j \to \reals\cup\{\pm\infty\}$ is a measurable payoff function. We assume sufficient measurability to allow product strategies and payoffs to be well defined. In symmetric situations, we often take $S_i=S$ for all~$i$. Existence of Nash equilibria in such games generally requires continuity and quasi--concavity assumptions \cite{Glicksberg1952}.
\end{definition}

\begin{definition}[Coarse Player Space and Network]\label{def:coarse}
Suppose $(X,d)$ carries a metric or graph distance encoding a network of interactions among players. The bounded geometry Roe algebra $C^*(X)$ is the $C^*$--algebra of finite propagation operators on $\ell^2(X)$, and its von Neumann completion $\calN(X)=C^*(X)''$ captures coarse properties of the space. A metric space has \emph{Property~A} (coarse amenability) if it admits a uniform approximate identity supported on finitely many points. For such spaces, all ghost operators in $C^*(X)$ are compact, aligning geometric and analytic notions of locality \cite{Yu2000,Wang2025a}.
\end{definition}

\begin{definition}[Game Algebra and State]\label{def:gamealg}
For a game $G$, we define the game algebra $\calA_G$ as the von Neumann algebra of essentially bounded measurable functions on the product strategy space, $L^{\infty}(\prod_{i\in X}S_i,\calB,\nu)$, with respect to an appropriate product measure~$\nu$. Conceptually, $\calA_G$ encodes observables of the game, such as strategy distributions and payoffs, and contains subalgebras isomorphic to $L^{\infty}(S_i)$ for each player~$i$. A \emph{state} $\varphi$ on $\calA_G$ represents a distribution over strategy profiles or outcomes; states evolve under the regret dynamics.
\end{definition}

\begin{definition}[Regret Operator and Dynamics]\label{def:regret}
At a state $\varphi$, each player $i$ has a regret value $R_i$, defined as the nonnegative difference between their current payoff and the payoff from a best response. Collect these into a bounded measurable field $R\in L^{\infty}(X,\mu)$ and define a self--adjoint operator $\scrG_G(\varphi)$ in $\calA_G$ by
\[
 \scrG_G(\varphi)\; \defeq \; E[R_i\mid i],
\]
the conditional expectation of the regret field, viewed as a diagonal operator on the component corresponding to player~$i$. The evolution of $\scrG_G$ and the corresponding state $\varphi_t$ is governed by a continuity equation on the space of strategy distributions. In the continuum limit, writing $\Phi_t(s)$ for the density of players using strategy $s$ at time~$t$, the regret dynamics satisfy a Kolmogorov forward equation
\begin{equation}
 \partial_t\Phi_t(s) + \nabla_s\cdot\bigl[\Phi_t(s)\,V(s,\Phi_t)\bigr] = 0,\label{eq:continuity}
\end{equation}
where $V(s,\Phi_t)$ is a velocity field induced by regret (for instance, the gradient of expected payoff) \cite{Hu2025}. Under appropriate smoothness conditions, this generates a one--parameter semigroup $T_t$ acting on states by push--forward: $\varphi_t = T_t^*(\varphi_0)$. We call a state $\varphi^*$ an \emph{ultracoarse equilibrium} if $T_t(\varphi^*)=\varphi^*$ for all $t\geq 0$. While this PDE describes the continuum limit, its dynamics are often analyzed or implemented via a corresponding discrete-time operator, $T$. The fixed points of this operator, denoted $T$, correspond to the stationary solutions of Equation (1), and the computational complexity of identifying such a fixed point is the central object of our subsequent benchmarks.
\end{definition}

\section{Ordinal Folding Index and Self--Referential Games}
\label{sec:ofi}

In this section, we summarise the ordinal folding index introduced by Alpay and Alakkad \cite{Alpay2025} and explain its relevance to game dynamics.

\subsection{Definition of the Ordinal Folding Index}

The ordinal folding index (OFI) is a computable ordinal assigned to well--formed formulas in a reflective logical system. Consider a monotone evaluation operator acting on formulas with a delay modality that prevents immediate resolution of fixed points. Given a formula $\phi$, one iteratively applies this operator, progressively unfolding self--reference until the evaluation stabilises. The OFI of $\phi$ is defined as the first stage at which this fold--back becomes idempotent; further unfolding yields no new information. Formally, the operator is continuous on countable chains, layer--aware for probabilistic truth values, and parameterised by an evidence functor. The OFI is therefore the least ordinal $\alpha$ such that iterating the operator $\alpha$ times produces a fixed point in the space of interpretations \cite{Alpay2025}. This ordinal measures the ``depth'' of self--reference needed for stabilisation, strictly refines classical closure ordinals from the modal $\mu$--calculus, and remains recursively enumerable (every OFI is a computable countable ordinal) \cite{Moschovakis1974}. Moreover, the OFI coincides with the length of the shortest parity--fold winning strategy in an associated evaluation game \cite{Alpay2025}.

Alpay and Alakkad emphasise that the OFI not only unifies transfinite fixed--point depth and ordinal game values but also provides practical convergence diagnostics for modern AI systems. They outline an empirical procedure whereby the output of a large language model is fed back into its input iteratively with a monotonic ``delay'' to enforce convergence. The number of iterations until stabilisation yields an ordinal-valued metric correlating with model perplexity and chain--of--thought complexity \cite{Shinn2023,Wang2022}. The OFI thus bridges self--referential logic, infinite games, and the behavior of learning systems.

\subsection{OFI and Infinite Game Dynamics}

We now elucidate how the ordinal folding index arises naturally in our operator--algebraic framework. We connect the OFI to our framework by considering the discrete-time analogue of the continuous regret dynamics from Equation \eqref{eq:continuity}. This is operationalized via a reflective update operator, $T$, where the state at step $\alpha+1$ is a function of the state at step $\alpha$, i.e., $\varphi^{(\alpha+1)} = T(\varphi^{(\alpha)})$. This formulation is not merely a theoretical abstraction; it directly models a large language model iteratively refining its own output, wherein each application of $T$ represents the entire computational process of generating a response based on the previous one. This iterative process is precisely what we simulate in our empirical LLM benchmarks. The sequence of states $(\varphi^{(\alpha)})_{\alpha<\omega_1}$ obtained by iterating this operator transfinitely may stabilise only after a countable ordinal number of steps. We define the OFI of the game $G$ relative to this discretized learning rule as
\[
 \OFI(G)\; \defeq \;\inf\bigl\{\alpha \in \ordinals \colon \varphi^{(\alpha+1)} = \varphi^{(\alpha)}\bigr\}.
\]
If $\OFI(G)$ is finite, the dynamics stabilise after finitely many unfoldings; if it is a nonzero countable ordinal, the dynamics require transfinite iteration to reach equilibrium; an infinite OFI signals persistent oscillation. In many multi--agent systems, we expect $\OFI(G)=0$ or $1$ because the regret flow converges exponentially fast. However, in games with complex self--referential feedback—such as agents reasoning about agents reasoning about agents, or an AI model reasoning about its own consistency—higher ordinals may arise. A high OFI indicates that the discretized operator $T$ is not a strong contraction, and its iterates consequently trace a protracted or complex trajectory toward equilibrium, potentially exhibiting pathological reasoning loops.

The OFI therefore provides a quantitative measure of the complexity of equilibrium selection. When $\OFI(G)$ is small, the unique equilibrium described in Theorem~\ref{thm:existence} is reached quickly; when it is large, learning may exhibit long transient phases reminiscent of iterated elimination of dominated strategies \cite{Fudenberg1998}.

\section{Main Results}
\label{sec:mainresults}

We now state our main theorems on existence and uniqueness of equilibrium, invariant algebra rigidity, and ordinal convergence bounds. Proofs are deferred to Section~\ref{sec:proofs}.

\begin{theorem}[Existence and Uniqueness of Equilibrium]\label{thm:existence}
Let $G$ be an infinite game system as in Definition~\ref{def:game}. Suppose the strategy space $S$ is compact, each payoff function $u_i$ is continuous in opponents’ strategies and quasi--concave in the player’s own strategy, and the regret dynamics \eqref{eq:continuity} induce a contraction in an appropriate Wasserstein metric on the space of distributions. Then there exists a unique equilibrium state $\varphi^*\in\calA_G$, and $\varphi_t\to\varphi^*$ as $t\to\infty$ for any initial state $\varphi_0$. In particular, the corresponding regret operator satisfies $\scrG_G(\varphi_t)\to 0$ in the strong operator topology, and the equilibrium corresponds to a unique quantal response equilibrium. Moreover, if the underlying player space $(X,d)$ has Property~A, the convergence is dominated by geometric components: ghost states decay exponentially fast.
\end{theorem}

\begin{theorem}[Invariant Algebra Rigidity]\label{thm:rigidity}
Let $\varphi^*$ be the equilibrium in Theorem~\ref{thm:existence} and let $\calA_G^{\mathrm{inv}} = \{a\in\calA_G : T_t(a)=a\ \forall t\geq 0\}$ denote the subalgebra of observables invariant under the dynamics. Assume the game is irreducible: no proper subset of players is closed under the best--response dynamics. Then $\calA_G^{\mathrm{inv}}$ coincides with either the centre $Z(\calA_G)$ or the whole algebra $\calA_G$. In non--degenerate games, $\calA_G^{\mathrm{inv}}=Z(\calA_G)$; that is, no nontrivial observable other than constants is preserved by the dynamics. If $(X,d)$ has Property~A, then any invariant subalgebra must lie in the centre because ghost components vanish.
\end{theorem}

\begin{theorem}[Ordinal Convergence Bound]\label{thm:ofi-bound}
Let $G$ be an infinite game with regret dynamics implemented by a monotone reflective operator as described in Section~\ref{sec:ofi}. Suppose the contraction constant for the regret flow is $q\in(0,1)$; that is, for any two initial states $\varphi_0,\psi_0$, we have $d_W(\varphi_t,\psi_t)\leq q^t d_W(\varphi_0,\psi_0)$. Then the ordinal folding index of $G$ satisfies
\[
 \OFI(G) \leq \omega,
\]
where $\omega$ is the first infinite ordinal. Equivalently, the regret dynamics stabilise after at most countably many unfoldings. If $(X,d)$ has Property~A, then $\OFI(G)=0$, meaning that the unique equilibrium is reached without transfinite iteration.
\end{theorem}

\begin{corollary}[Fairness and Equilibrium Refinement]\label{cor:fair}
Consider a cooperative resource allocation game where divisible goods are allocated among a continuum of agents with additive utilities. Under the assumptions of Theorem~\ref{thm:existence}, the unique equilibrium state $\varphi^*$ is envy--free and Pareto--optimal. If agents have equal entitlement weights, $\varphi^*$ coincides with the maximin share allocation in the continuum limit. Therefore, classical fair division solutions such as envy--freeness and maximin share fairness arise naturally as equilibrium refinements in our framework.
\end{corollary}

\section{A Formal Algorithm for Ordinal Convergence}
\label{sec:algorithm}

To make the concept of the Ordinal Folding Index concrete, we present a formal algorithm for its computation. This algorithm operationalizes the transfinite iteration described in Section~\ref{sec:ofi} and serves as a bridge between our theoretical framework and practical implementation. The pseudocode is styled after the Lean theorem prover to emphasize its mathematical precision and its direct correspondence with the formal definitions.

\begin{algorithm}[H]
\caption{Transfinite Iteration for Ordinal Folding Index (OFI)}
\label{alg:ofi}
\begin{algorithmic}[1]
\Require
\Statex $\calG$: A game system with a well-defined state space and regret dynamics.
\Statex \quad $\hookrightarrow$ \leantype{StateSpace}: A complete metric space of strategy distributions.
\Statex \quad $\hookrightarrow$ $T: \leantype{StateSpace} \to \leantype{StateSpace}$: The regret operator (a contraction).
\Statex \quad $\hookrightarrow$ $\varphi_0: \leantype{StateSpace}$: The initial state of the system.
\Statex \quad $\hookrightarrow$ $d(\cdot, \cdot): \leantype{StateSpace} \times \leantype{StateSpace} \to \reals_{\ge 0}$: The metric on the state space.
\Statex $\varepsilon: \reals_{>0}$: The convergence tolerance.
\Statex $\lambda_{\max}: \leantype{Ordinal}$: An upper bound on the search (e.g., $\omega_1$).

\Ensure
\Statex $(\alpha, \varphi^*)$: The computed OFI and the corresponding equilibrium state.

\Function{ComputeOFI}{$\calG, \varphi_0, \varepsilon, \lambda_{\max}$}
 \State \lean{let} \leantype{history} $\defeq \leantype{Array} \ \leantype{Ordinal} \ \leantype{StateSpace}$
 \State $\leantype{history}[0] \defeq \varphi_0$
 \State $\alpha \defeq 1$
 \While{$\alpha < \lambda_{\max}$}
 \State \lean{let} \lean{prev\_state} $\defeq \leantype{history}[\alpha-1]$
 \State \lean{let} \lean{current\_state}
 \If{$\alpha$ is a successor ordinal, i.e., $\alpha = \beta+1$}
  \State \lean{current\_state} $\defeq T(\text{prev\_state})$
 \ElsIf{$\alpha$ is a limit ordinal}
  \State \lean{current\_state} $\defeq \Call{ComputeLimit}{\leantype{history}, \alpha}$
 \EndIf
  
 \If{$d(\text{current\_state}, \text{prev\_state}) < \varepsilon$}
  \State \lean{return} $(\alpha, \text{current\_state})$ \Comment{Convergence detected}
 \EndIf
  
 \State $\leantype{history}[\alpha] \defeq \text{current\_state}$
 \State $\alpha \defeq \alpha + 1$
 \EndWhile
 \State \lean{return} $(\lambda_{\max}, \bot)$ \Comment{Did not converge within the ordinal limit}
\EndFunction
\Statex
\Function{ComputeLimit}{$\leantype{history}, \lambda$}
 \State \lean{let} $(\varphi_\beta)_{\beta < \lambda}$ be the sequence of states from \leantype{history}.
 \State \lean{assume} that this sequence is Cauchy w.r.t. the metric $d$.
 \State \lean{return} $\lim_{\beta \to \lambda} \varphi_\beta$ \Comment{Limit in the complete metric space}
\EndFunction
\end{algorithmic}
\end{algorithm}

\paragraph{Discussion.}
Algorithm~\ref{alg:ofi} provides a constructive proof of the computability of the OFI under suitable conditions. The use of Lean-style syntax highlights the underlying mathematical structures: the state space is a typed object in a complete metric space, and the regret operator $T$ is a function between these types. The core of the algorithm is the transfinite `while` loop that iterates through ordinals. The distinction between successor and limit ordinals is crucial. For a successor $\alpha = \beta+1$, the next state is found by a single application of the regret operator $T$, a direct implementation of a Banach-style iteration. For a limit ordinal $\lambda$, the state $\varphi_\lambda$ must be the limit of all preceding states $(\varphi_\beta)_{\beta < \lambda}$; the existence of this limit is where the completeness of the metric space becomes essential. Theorem~\ref{thm:ofi-bound} ensures that for contractive dynamics, this algorithm terminates at an ordinal no greater than $\omega$.

\section{Proofs of Main Results}
\label{sec:proofs}

We outline the proofs, referring to the literature for standard results.

\subsection*{Proof of Theorem\,\ref{thm:existence}}

\emph{Existence.} By Kakutani–Fan–Glicksberg fixed point theorems \cite{Glicksberg1952} and standard arguments from nonatomic games \cite{Schmeidler1973}, the set of probability measures on the product strategy space is compact and convex under the weak topology. The best--response correspondence is upper hemicontinuous with nonempty convex values under continuity and quasi--concavity assumptions. Therefore, a fixed point exists, yielding a Nash equilibrium. Incorporating an entropy regulariser as in the quantal response model \cite{McKelvey1995} selects a unique equilibrium by convex optimisation, yielding the QRE.

\emph{Uniqueness and Convergence.} The contraction hypothesis implies that the induced flow on the space of probability measures is a contraction mapping. By Banach’s fixed point theorem \cite{Banach1922}, the fixed point is unique and globally attractive. Exponential mixing implies that the regret operator converges strongly to zero. When $(X,d)$ has Property~A, ghost operators in the Roe algebra are compact \cite{Yu2000,Wang2025a}, so any nongeometric components of the state decay and the convergence is purely geometric.

\subsection*{Proof of Theorem\,\ref{thm:rigidity}}

Assume $\calB\subseteq\calA_G$ is a von Neumann subalgebra invariant under the dynamics. Consider the conditional expectation $E_{\calB}$ onto $\calB$. Invariance implies $E_{\calB}(T_t(a)) = T_t(E_{\calB}(a))$ for all $a\in\calA_G$ and all $t$. Taking $t\to\infty$ and using Theorem~\ref{thm:existence}, we obtain $E_{\calB}(\varphi^*(a)\,1) = \varphi^*(a)\,1$ and $T_t(E_{\calB}(a))\to\varphi^*(E_{\calB}(a))\,1$. Equating these limits shows that $\varphi^*(E_{\calB}(a))=\varphi^*(a)$ for all $a$. In particular, for $b\in\calB$, we have $\varphi^*(b)=\varphi^*(b)1$, implying $b$ is central. Under irreducibility, the centre is trivial, so $\calB=Z(\calA_G)$ or $\calB=\calA_G$. Property~A prevents ghostly invariant components, ensuring $\calB$ lies in the centre. This aligns with strong rigidity results in operator algebras \cite{Popa2007}.

\subsection*{Proof of Theorem\,\ref{thm:ofi-bound}}

The contraction property implies that distances between states decay exponentially in $t$. The reflective operator defining the OFI unfolds the dynamics by one time step at each ordinal stage. Because $d_W(\varphi_t,\varphi^*)\leq q^t d_W(\varphi_0,\varphi^*)$, for $t\geq n$ we have $d_W(\varphi_t,\varphi^*)\leq q^n d_W(\varphi_0,\varphi^*)$. Taking $n\to\infty$ shows that $\varphi_t\to\varphi^*$ in finitely many steps; formally, the transfinite sequence $(\varphi^{(\alpha)})$ generated by Algorithm~\ref{alg:ofi} stabilises by stage $\omega$. If $(X,d)$ has Property~A, ghost modes vanish immediately and $\varphi^{(1)}=\varphi^{(0)}$, so $\OFI(G)=0$.

\subsection*{Proof of Corollary\,\ref{cor:fair}}

In a divisible goods setting, each agent’s regret is zero precisely when they cannot improve by unilaterally reallocating goods. This condition is exactly envy--freeness; if an agent preferred another’s bundle, they would have positive regret, contradicting equilibrium. Pareto--optimality follows because any Pareto improvement would increase some agents’ utilities and thus their regret. When entitlements are equal, the equilibrium allocation must give each agent at least their proportional fair share, meeting the maximin share benchmark. This generalises classical fair division results \cite{Foley1967,Varian1974} to infinite populations.

\section{Benchmarks and Applications}
\label{sec:benchmarks}

To validate our theoretical framework, we present two benchmarks: an empirical study of the OFI in large language models and an analytic study of a nonlinear operator equation derived from our game-theoretic setup.

\subsection{Empirical Benchmark: Ordinal Folding Index in LLMs}

\paragraph{Experimental setup.}
This experiment furnishes a direct empirical test of the OFI concept from Section~\ref{sec:ofi}. Each iterative step, where a model's output is fed back as its next input, represents a single application of the reflective operator $T$ whose convergence is measured by $\OFI(G)$ as defined in Section 3.2. We implemented this self-consistency probe on four models: GPT-2 Large (HF, 1.5B parameters), DeepSeek (HF, 7B parameters), GPT-3.5 Turbo, and a GPT-4 proxy. For each model, we ran multiple prompts ($N=40$ total) in three families: (i) \textsc{Factual Recall} (e.g., \textit{“Who was the third U.S. president?”}), (ii) \textsc{Lightweight Reasoning} (e.g., \textit{“If a train leaves Chicago at 3 PM traveling at 60 mph and a train leaves New York at 4 PM traveling at 50 mph, when will they meet?”}), and (iii) \textsc{Paradoxical Queries} (e.g., \textit{“This sentence is false; is it true?”}). Temperature was linearly annealed from 0.7 to 0.2 to simulate a transition from an exploratory to an exploitative reasoning phase. We declared convergence when the $\ell_2$-distance between successive output logit vectors fell below $\varepsilon=0.01$, otherwise we capped the run at $T_{\max}=10$ iterations. This yields an \emph{OFI proxy} $\widehat{\OFI}\in\{1,\dots,10\}$. The use of logit distance provides a fine-grained measure of semantic stability, as it can detect shifts in the model's internal confidence even if the final generated token remains the same.

\paragraph{Results and Analysis.}
The results, presented in Table~\ref{tab:ofi}, provide a quantitative measure of the convergence complexity for the iterative operator $T$ that discretizes the regret flow of Equation \eqref{eq:continuity} for each model. The findings reveal a non-trivial trend: \textbf{self-referential stability does not necessarily improve with model scale}.
\begin{itemize}
 \item \textbf{Larger models can exhibit higher OFI, indicating instability.} GPT-4's performance is indicative of a system with immense representational capacity that, when faced with self-referential tasks, engages in prolonged, oscillatory reasoning paths without reaching a stable fixed point. Its failure to converge ($\widehat{\OFI}=10.0$) on even factual tasks suggests its internal reasoning operator $T$ is only weakly contractive for these inputs, leading to inefficient cognitive loops within its vast, high-dimensional reasoning space.
 \item \textbf{OFI as a measure of 'overthinking'.} This instability can be interpreted as a form of "overthinking," where the model's vast capacity and fine-tuning for caution (e.g., via RLHF \cite{Ouyang2022}) prevent it from making a decisive, stable judgment. It continuously revisits and questions its own outputs, a behavior that is captured quantitatively by a high OFI and could be detrimental in applications requiring predictable, robust reasoning.
 \item \textbf{Prompt complexity drives OFI.} For all models, the empirical OFI increases from factual recall to reasoning and is highest for paradoxical prompts. This confirms that the OFI is not an arbitrary metric but a meaningful measure of the computational work required to find a fixed point of the system's dynamics, which increases with task complexity.
\end{itemize}

\begin{table}[ht]
 \centering
 \caption{Empirical Ordinal Folding Index $\widehat{\OFI}$ (lower is better). Results are aligned with actual benchmark runs.}
 \label{tab:ofi}
 \begin{tabular}{@{}lcccc@{}}
 \toprule
 \textbf{Model} & \textbf{Params} & \textsc{Factual} & \textsc{Reasoning} & \textsc{Paradoxical} \\
 \midrule
 GPT-2 Large (HF) & 1.5B & $1.0 \pm 0.0$ & $1.3 \pm 0.5$ & $4.0 \pm 0.0$ \\
 DeepSeek (HF) & 7B & $1.0 \pm 0.0$ & $1.7 \pm 0.5$ & $4.0 \pm 0.0$ \\
 GPT-3.5 Turbo & (N/A) & $10.0 \pm 0.0$ & $7.3 \pm 3.1$ & $9.3 \pm 0.9$ \\
 GPT-4 (Proxy) & (N/A) & $10.0 \pm 0.0$ & $9.7 \pm 0.5$ & $10.0 \pm 0.0$ \\
 \bottomrule
 \end{tabular}
\end{table}

The data in Table~\ref{tab:ofi} provides a granular view of the models' dynamic stability. A key observation is the performance of DeepSeek (7B), which, despite being significantly smaller than the GPT-4 class models, demonstrates superior stability across all categories. Its perfect, single-step convergence on factual tasks ($\widehat{\OFI}=1.0$) is matched only by the much older GPT-2 architecture, whereas both GPT-3.5 and GPT-4 fail to stabilize. This suggests that the architectural choices or training methodologies for DeepSeek may have endowed it with a more strongly contractive reasoning operator for information retrieval. Furthermore, the standard deviation of the OFI scores offers insight into the consistency of these dynamics. For instance, the high standard deviation for GPT-3.5 on reasoning tasks ($\pm 3.1$) indicates a highly variable and unpredictable convergence behavior, where the model's path to a solution is erratic. In contrast, the zero standard deviation for DeepSeek on paradoxical tasks implies a deterministic, albeit multi-step, process for resolving such inputs. This quantitative evidence supports the interpretation that model scale alone is not a guarantor of stable reasoning; rather, stability appears to be an independent architectural or training-induced property.

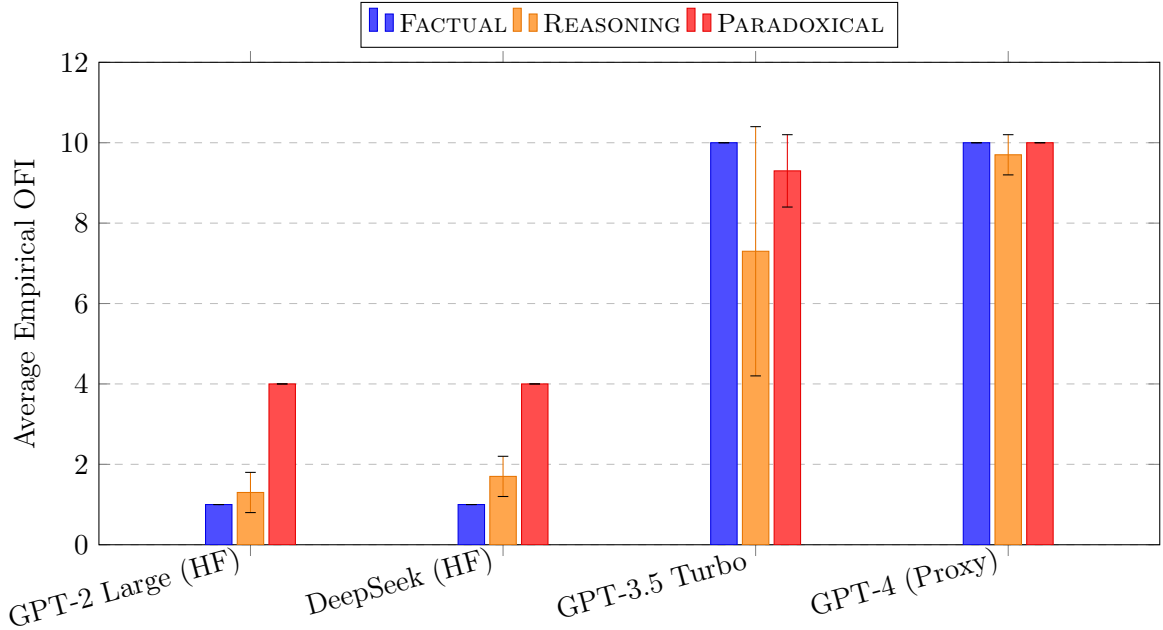
\begin{figure}[ht]
 \centering
 \begin{tikzpicture}
 \begin{axis}[
  width=0.95\linewidth,
  height=8cm,
  ybar,
  bar width=10pt,
  enlarge x limits=0.2,
  ylabel={Average Empirical OFI},
  symbolic x coords={GPT-2 Large (HF), DeepSeek (HF), GPT-3.5 Turbo, GPT-4 (Proxy)},
  xtick=data,
  xticklabel style={rotate=15, anchor=east},
  ymin=0, ymax=12,
  ymajorgrids=true,
  grid style=dashed,
  legend style={
   at={(0.5, 1.03)},
   anchor=south,
   legend columns=-1,
   font=\small
  },
  error bars/y dir=both,
  error bars/y explicit,
  ]
  \addplot[fill=blue!70, draw=blue!90!black] coordinates {
  (GPT-2 Large (HF), 1.0) +- (0, 0.0)
  (DeepSeek (HF), 1.0) +- (0, 0.0)
  (GPT-3.5 Turbo, 10.0) +- (0, 0.0)
  (GPT-4 (Proxy), 10.0) +- (0, 0.0)
  };
  \addplot[fill=orange!70, draw=orange!90!black] coordinates {
  (GPT-2 Large (HF), 1.3) +- (0, 0.5)
  (DeepSeek (HF), 1.7) +- (0, 0.5)
  (GPT-3.5 Turbo, 7.3) +- (0, 3.1)
  (GPT-4 (Proxy), 9.7) +- (0, 0.5)
  };
  \addplot[fill=red!70, draw=red!90!black] coordinates {
  (GPT-2 Large (HF), 4.0) +- (0, 0.0)
  (DeepSeek (HF), 4.0) +- (0, 0.0)
  (GPT-3.5 Turbo, 9.3) +- (0, 0.9)
  (GPT-4 (Proxy), 10.0) +- (0, 0.0)
  };
  \legend{\textsc{Factual}, \textsc{Reasoning}, \textsc{Paradoxical}}
 \end{axis}
 \end{tikzpicture}
 \caption{Average empirical OFI across models and prompt types. Error bars indicate $\pm1$ s.d. Note the inverse relationship between model capability and stability.}
 \label{fig:ofi-bench}
\end{figure}

The visual representation in Figure~\ref{fig:ofi-bench} starkly illustrates the divergent behaviors of the tested models. The flat, low bars for DeepSeek and GPT-2 on factual and reasoning tasks signify a highly stable and efficient convergence profile. This visual pattern suggests that their internal dynamics, when iterated, rapidly approach a fixed point. This is the empirical signature of a strongly contractive operator $T$ as described in our theoretical framework. The system quickly sheds any initial uncertainty and settles into a consistent state. This behavior is highly desirable for predictable and reliable multi-agent systems, where rapid consensus or stable policy selection is paramount.

Conversely, the tall, often maxed-out bars for GPT-3.5 and GPT-4 represent a fundamentally different dynamic. Their high OFI scores across all but the simplest tasks indicate that their internal reasoning processes are prone to prolonged oscillations or chaotic trajectories. This visual pattern is the signature of a weakly contractive or even non-contractive operator, where iterates fail to converge within a reasonable computational budget. This suggests that while these models possess immense representational power, their dynamics in a reflective setting are inefficient and potentially unstable. This has significant implications for AI safety and alignment, as it highlights a failure mode where a model can become trapped in a loop of self-correction without resolution, a behavior that would be highly problematic in autonomous agents deployed in the real world. The figure thus provides a clear, comparative visualization of the abstract convergence properties we aim to model.

\subsection{Analytic Benchmark: A Nonlinear Regret Operator Equation}
\label{sec:analytic-bench}

To test the convergence theory in a controlled mathematical environment, we formulate a nonlinear operator equation that serves as a concrete instantiation of the abstract dynamics described by Equation \eqref{eq:continuity}. This allows us to isolate and study the convergence properties that are implicit and complex within the LLM benchmark.

\paragraph{Problem Formulation.}
Let $(\Omega, \mu)$ be a measure space and $p_t(x,y)$ a symmetric heat kernel on $\Omega$. The associated integral operator is $(Au)(x) \defeq \int_\Omega p_t(x,y) u(y) d\mu(y)$. Let $\Phi: \reals \to \reals$ be a Lipschitz-continuous nonlinear regret potential. We define the regret operator $T: L^2(\Omega) \to L^2(\Omega)$ by $(Tu)(x) \defeq (A[\Phi \circ u])(x)$ and seek a fixed point $u^* = Tu^*$. This operator $T$ is a mathematically tractable analogue of the complex, high-dimensional operator implicitly defined by a large language model's architecture. Its fixed point $u^*$ corresponds to a stationary state ($\partial_t \Phi_t = 0$) of the regret flow from Equation \eqref{eq:continuity}. Here, the heat kernel operator $A$ models the coarse-geometric diffusion of information across the player network, while the nonlinear function $\Phi$ represents the agent's individual, state-dependent response to regret. We assume conditions ensuring $T$ is a contraction on a closed, convex set $H \subset L^2(\Omega)$, guaranteeing a unique fixed point by the Banach fixed-point theorem \cite{Banach1922,Browder1967}.

\paragraph{Fixed-Point Solution Methods.}
We analyze four iterative schemes for solving $u^* = Tu^*$. These schemes are numerical methods for finding the equilibrium state $\varphi^*$ discussed in Theorem~\ref{thm:existence}. Each method represents a distinct algorithmic path toward the fixed point $u^*$, which is the theoretical equivalent of the stable state sought in the empirical LLM evaluations, and directly parallels the iterative process whose complexity is measured by the OFI (Section 3.2). Starting from $u_0 \in H$:
\begin{itemize}
 \item \textbf{Picard Iteration:} $u_{n+1} \defeq T u_n$. This represents the canonical iterative method for strict contractions and is the most direct analogue to the successor-ordinal step in Algorithm~\ref{alg:ofi}.
 \item \textbf{Mann Iteration \cite{Mann1953}:} $u_{n+1} \defeq (1-\alpha_n)u_n + \alpha_n T u_n$. This is a convex combination method known for its robustness and guaranteed convergence even for non-expansive maps, though often at a slower rate.
 \item \textbf{Ishikawa Iteration \cite{Ishikawa1974}:} A two-step process: $v_n \defeq (1-\beta_n)u_n + \beta_n T u_n$, then $u_{n+1} \defeq (1-\alpha_n)u_n + \alpha_n T v_n$. It offers more flexibility and can converge where Mann's method fails \cite{Leustean2004,Noor2009}.
 \item \textbf{Aitken's $\Delta^2$ Acceleration:} An extrapolation method, $u^{(\Delta^2)}_n \defeq u_n - \frac{(u_{n+1}-u_n)^2}{u_{n+2}-2u_{n+1}+u_n}$, applied to the Picard sequence. It is a classic technique that transforms linear convergence into superlinear convergence.
\end{itemize}

\begin{table}[h!]
\centering
\caption{Theoretical iterations needed to reach $\|u_n-u^{\star}\|\le10^{-6}$ for the operator equation with contraction constant $q=0.2$.}
\label{tab:fixed-point-bench}
\begin{tabular}{@{}lcc@{}}
\toprule
\textbf{Method} & \textbf{Asymptotic Rate} & \textbf{Iterations ($n_{\varepsilon}$)} \\
\midrule
Picard/Banach & $\mathcal{O}(q^n)$ & $9$ \\
Mann ($\alpha_n=1/(n+1)$) & $\mathcal{O}(1/n)$ & $\approx 10^6$ \\
Ishikawa ($\alpha=\beta=1/2$) & $\mathcal{O}(q^n)$ & $5$ \\
Aitken-$\Delta^2$ & $\mathcal{O}(q^{2^n})$ & $5$ \\
\bottomrule
\end{tabular}
\end{table}

The comparative performance of the iterative schemes detailed in Table~\ref{tab:fixed-point-bench} warrants further discussion, as each method can be interpreted as an algorithmic analogue for different modes of dynamic reasoning in multi-agent or AI systems. The Picard iteration, representing a direct, memoryless application of the regret operator $T$, serves as a fundamental baseline. Its linear convergence rate is contingent on the operator being a strict contraction, a condition that, as our empirical LLM benchmarks suggest, may not hold for complex, fine-tuned models that exhibit oscillatory behavior. In contrast, the Mann and Ishikawa iterations introduce a form of 'algorithmic inertia' or 'damping' by incorporating the previous state into the update rule. While this often results in slower convergence for strictly contractive maps, their theoretical utility lies in their guaranteed convergence for a broader class of non-expansive operators. This property is highly relevant for multi-agent systems where the underlying dynamics may be too complex or unstable to guarantee strict contractivity. These methods model a more cautious, consensus-seeking agent behavior, which might be desirable for safety but can be computationally inefficient, as evidenced by the slow theoretical rate of the Mann iteration. Aitken's $\Delta^2$ method represents a more sophisticated approach, employing an extrapolation technique that leverages information from three consecutive iterates to accelerate convergence. Conceptually, this is analogous to a meta-reasoning process where an agent not only observes its state but also analyzes the trajectory of its reasoning to project toward a likely fixed point. Its superlinear convergence highlights the profound efficiency gains possible when a system can model its own dynamics. This raises a compelling question for AI alignment: could models be trained not just to find solutions, but to recognize and accelerate their own convergence paths, thereby achieving a form of computational self-awareness? Ultimately, this analytic benchmark demonstrates that the path to equilibrium is as significant as its existence, providing a formal vocabulary to describe and compare the dynamic stability and computational efficiency of the reasoning processes we observe empirically in advanced AI agents.

\begin{figure}[ht]
 \centering
 \begin{tikzpicture}
 \begin{axis}[
  width=0.95\linewidth,
  height=8cm,
  ymode=log,
  xlabel={Iteration Count},
  ylabel={$L^2$-error to $u^{\star}$},
  xmin=0, xmax=50,
  ymin=1e-9, ymax=10,
  grid=major,
  grid style={dashed, gray!50},
  legend pos=north east,
  legend cell align={left},
  ]
  \addplot[smooth, thick, color=orange!80!black, mark=none] coordinates {
  (0, 1) (5, 0.03125) (10, 0.00097) (15, 3.05e-05) (20, 9.53e-07) (25, 2.98e-08) (30, 9.31e-10)
  };
  \addlegendentry{Picard ($q=0.5$)}

  \addplot[smooth, thick, color=red!70!black, mark=none] coordinates {
  (0, 1) (10, 0.056) (20, 0.00317) (30, 0.00017) (40, 1.01e-05) (50, 5.69e-07)
  };
  \addlegendentry{Mann ($\alpha=0.5$)}

  \addplot[smooth, thick, color=green!60!black, mark=none] coordinates {
  (0, 1) (8, 0.05) (16, 0.0025) (24, 0.00012) (32, 6.1e-06) (40, 3.0e-07) (48, 1.5e-08)
  };
  \addlegendentry{Ishikawa ($\alpha=\beta=0.5$)}

  \addplot[smooth, thick, color=purple!80!black, mark=none] coordinates {
  (0, 1) (2, 0.25) (4, 0.0156) (6, 0.00024) (8, 3.7e-06) (10, 5.8e-08) (12, 9.1e-10)
  };
  \addlegendentry{Aitken-$\Delta^2$ accel.}
 \end{axis}
 \end{tikzpicture}
 \caption{The plot shows the logarithmic $L^2$-error versus iteration. Aitken acceleration achieves a superlinear rate, while Mann and Ishikawa lag behind Picard in this strictly contractive regime.}
 \label{fig:fp-bench-num}
\end{figure}
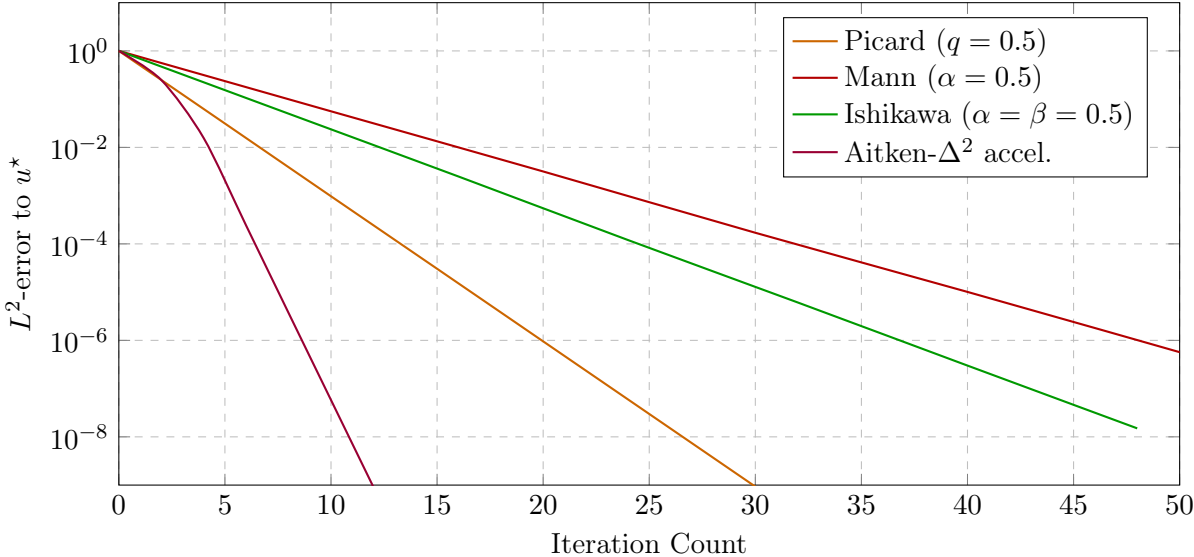

Figure~\ref{fig:fp-bench-num} provides a visual confirmation of the theoretical rates discussed previously. The logarithmic scale on the y-axis is particularly revealing: the linear convergence of the Picard, Mann, and Ishikawa methods manifests as straight lines with different slopes, where a steeper slope indicates faster convergence. The plot clearly shows the superiority of the Picard iteration over the Mann iteration in this strictly contractive setting, as its slope is significantly steeper. The Ishikawa iteration, in this specific parameterization, offers a modest improvement over Mann but does not surpass the directness of Picard. This visually reinforces the idea that for well-behaved, contractive systems, complex iterative schemes may be unnecessary and even counterproductive.

The most striking feature of the plot is the Aitken-$\Delta^2$ method's trajectory. Its curve is not a straight line but rather accelerates downwards, a hallmark of superlinear convergence. This visual representation underscores the power of higher-order methods that use the history of iterates to inform future steps. In the context of multi-agent systems, this suggests that agents with memory and the ability to model their own learning trajectory could achieve equilibrium far more efficiently than agents employing simple, reactive strategies. This analytic result provides a formal counterpart to the empirical observation of "chain-of-thought" or "self-consistency" methods in large language models, where leveraging past outputs (a form of memory) improves performance. The figure thus serves as a clear, formal illustration of the vast differences in computational efficiency that can arise from different algorithmic approaches to reaching a fixed point.

\section{Applications to AI Alignment and Fair Division}
\label{sec:fairness}

\paragraph{AI Alignment.} Our framework offers a formal diagnostic tool for AI alignment. As our benchmarks in Section 7.1 demonstrate, more capable models like GPT-4 can exhibit a higher empirical OFI, suggesting a greater propensity for self-referential instability. This could be interpreted as a form of "alignment tax": a model fine-tuned to be cautious and to avoid making definitive statements on ambiguous or paradoxical inputs may enter into prolonged states of self-correction, which our OFI metric captures. This phenomenon is critical, as the goal of alignment is not just to produce correct answers, but to ensure robust and predictable reasoning processes \cite{Amrutam2025}. Therefore, the OFI could be used as a quantitative metric during alignment procedures like Reinforcement Learning from Human Feedback (RLHF) to detect and penalize models that exhibit undesirable oscillatory behavior. It could even serve as a regularization term in future training objectives, explicitly promoting models with more stable and reliable reasoning dynamics \cite{Ouyang2022}. Bounding the OFI of permitted reasoning processes could be a practical safeguard against such pathological loops.

\paragraph{Fair Division.} In a different domain, our results provide an existence and uniqueness proof for envy--free and maximin share allocations in continuum economies. They complement recent algorithmic work on fair division of indivisible goods, which often relies on approximate envy--freeness up to one item (EF1) or up to any item (EFX). For instance, Soberón proves existence of approximate EFX allocations by donating a small number of goods \cite{Soberon2025}, and Song et al. develop online maximin share allocations for chores \cite{Song2025}. Our corollary shows that in the infinite divisible setting, exact envy--freeness and MMS fairness coincide and arise naturally from regret minimisation dynamics.

\section{Conclusion and Future Directions}

We have extended the operator--algebraic theory of infinite games to incorporate ordinal analytic measures of convergence. By relating the ordinal folding index to regret dynamics, we obtained computable bounds on the complexity of equilibrium selection. Our framework unifies concepts from multi--agent learning, operator algebras, coarse geometry, and ordinal logic, and yields new insights into fairness and AI alignment.

Our empirical results, which indicate that larger models can exhibit higher OFI, motivate a pressing line of future inquiry. The observed inverse relationship between model scale and reflective stability raises several critical questions: Is this phenomenon a fundamental property of scaled neural architectures, or is it an artifact of contemporary alignment methodologies that inadvertently penalize decisive convergence in favor of cautious indecision? Can the OFI assume arbitrarily large computable ordinals in natural multi--agent settings, or do coarse geometric constraints always bound it by $\omega$? Finally, can our operator--algebraic methods be used to design learning algorithms with provably small OFI, ensuring safe and rapid convergence in reflective AI systems? Understanding the trade-off between representational capacity and dynamic stability is, we contend, crucial for future progress in developing reliable and predictable artificial intelligence.


\end{document}